\def\b#1{{\bf #1}}
\def\i#1{{\it #1}}

\documentstyle[11pt]{article}
\begin{document}

\title{On the Average of Triangular Numbers}
\author{Mario Catalani\\
Department of Economics, University of Torino\\
Via Po 53, 10124 Torino, Italy\\
mario.catalani@unito.it}
\date{}
\maketitle
\begin{abstract}
The problem we are dealing with is the following: find two sequences
$a_n$ and $b_n$
such that the average of the first $b_n$ triangular numbers (starting with
the triangular number 1) is still a triangular number, precisely
the $a_n$-th triangular number. We get also some side
results: for instance one of the sequence instrumental to finding
the asked for sequences turns out to be a bisection of
the sequence of the numerators
of continued fraction convergents to $\sqrt{3}$.
\end{abstract}

\bigskip
\bigskip
\bigskip

\noindent
The present note has been suggested by a problem proposed in the "Student
Problems" section of \i{The Mathematical Gazette} (see \cite{studentproblem}).
The problem we tackle is the following: find two sequences $a_n$ and $b_n$
such that the average of the first $b_n$ triangular numbers (starting with
the triangular number 1) is still a triangular number, precisely
the $a_n$-th triangular number.

\noindent
If we want that the average of the first $s$ triangular number be a
triangular number it has to hold
$${1\over s}\sum_{k=1}^s{k(k+1)\over 2}={r(r+1)\over 2},$$
for some positive integer $r$. This becomes
$${(s+1)(2s+4)\over 12}={r(r+1)\over 2},$$ that is, clearing the
fractions,
\begin{equation}
\label{eq:equazione}
s^2+3s+2=3r^2+3r.
\end{equation}
Solving for $s$ and considering that $s$ has to be a positive quantity
we have the solution
$$s={-3+\sqrt{1+12r+12r^2}\over 2}.$$
Now because $s$ has to be an integer we need $r$ such that $1+12r+12r^2$
is a perfect square; also $\sqrt{1+12r+12r^2}-3$ has to be even.

\noindent
Consider the nonhomogeneous recurrence
$w_n=w(k,\,r,\,s)$ defined by
$$w_n=4w_{n-1}-w_{n-2}+k,\quad w_0=r,\,w_1=s.$$
Let $\alpha$ and $\beta$ be the zeros of the polynomial $x^2-4x+1$,
that is $\alpha=2+\sqrt{3},\,\beta=2-\sqrt{3}$. Note that $\alpha+\beta=4,\,
\alpha\beta=1,\,\alpha-\beta=2\sqrt{3},\,\alpha^2+\beta^2=14$.

\noindent
It is useful to introduce the homogeneous sequence $L_n=w(0,\,2,4)$
given by
$$L_n=4L_{n-1}-L_{n-2}.$$
This sequence is the analogous in this context of the Lucas numbers.
The closed form of $L_n$ is
$$L_n=\alpha^n+\beta^n.$$
The generating function of $w_n$ is given by
$$g(x)={r+(s-5r)x+(k+4r-s)x^2\over (1-x)(1-4x+x^2)},$$
from which we get the closed form
\begin{equation}
w_n=-{k\over 2}+{s\alpha+(k+4r-s)\beta+k-2r\over 12}\alpha^n+
{s\beta+(k+4r-s)\alpha+k-2r\over 12}\beta^n.
\end{equation}
It may be useful to write the RHS in terms of the $L_n$ sequence
\begin{equation}
\label{eq:closed}
w_n=-{k\over 2}+{4s+k-2r\over 12}L_n+{k+4r-2s\over 12}L_{n-1}.
\end{equation}
Now introduce the sequence $a_n=w(1,\,0,\,1)$. $a_n$ is sequence A061278 in
\cite{sloane}: the first values are $0,\,1,\,5,\,20,\,76,\,285,\,\ldots$.
Using
Equation~\ref{eq:closed} we get
\begin{eqnarray*}
a_n&=&{5\over 12}L_n-{1\over 12}L_{n-1}-{1\over 2}\\
&=&{L_n+4L_n-L_{n-1}-6\over 12}\\
&=&{L_n+L_{n+1}-6\over 12}.
\end{eqnarray*}
Note that
\begin{eqnarray*}
L_n^2&=&(\alpha^n+\beta^n)^2\\
&=&\alpha^{2n}+\beta^{2n}+2\\
&=&L_{2n}+2,
\end{eqnarray*}
\begin{eqnarray*}
L_nL_{n+1}&=&(\alpha^n+\beta^n)(\alpha^{n+1}+\beta^{n+1})\\
&=&\alpha^{2n+1}+\alpha^n\beta^{n+1}+\alpha^{n+1}\beta^n+\beta^{2n+1}\\
&=&L_{2n+1}+\alpha+\beta\\
&=&L_{2n+1}+4.
\end{eqnarray*}
Then
\begin{eqnarray*}
1+12a_n+12a_n^2&=&{-12+L_{2n}+2L_{2n+1}+L_{2n+2}\over 12}\\
&=&{-12+6L_{2n+1}\over 12}\\
&=&-1+{L_{2n+1}\over 2},
\end{eqnarray*}
where we used the recurrence defining $L_{n}$. Now we are going to show
that this expression is a perfect square. Write
\begin{eqnarray*}
-1+{L_{2n+1}\over 2}&=&{1\over 2}\alpha^{2n+1}+{1\over 2}\beta^{2n+1}-1\\
&=&\left ({1\over \sqrt{2}}\alpha^{n+{1\over 2}}-
{1\over \sqrt{2}}\beta^{n+{1\over 2}}\right )^2.
\end{eqnarray*}
Hence it remains to prove that
$$u_n=
{1\over \sqrt{2}}\alpha^{n+{1\over 2}}-
{1\over \sqrt{2}}\beta^{n+{1\over 2}}$$
is a positive integer. This will be done by strong induction. First of all
note that
\begin{eqnarray*}
\left (\alpha^{1\over 2}-\beta^{1\over 2}\right )^2&=&
\alpha+\beta-2\\
&=&2,
\end{eqnarray*}
so that we will take $\alpha^{1\over 2}-\beta^{1\over 2}=\sqrt{2}$. Then
$u_0=1$. Next assume that $u_n$ is a positive integer for $n\le n_0$. It
follows that $u_{n_0}(\alpha+\beta)=4u_{n_0}$ is a positive integer. But
\begin{eqnarray*}
u_{n_0}(\alpha+\beta)&=&
{1\over \sqrt{2}}\left (\alpha^{n_0+{1\over 2}}-
\beta^{n_0+{1\over 2}}\right )(\alpha+\beta)\\
&=&{1\over \sqrt{2}}\alpha^{n_0+1+{1\over 2}}
+{1\over \sqrt{2}}\alpha^{n_0+{1\over 2}}\beta
-{1\over \sqrt{2}}\alpha\beta^{n_0+{1\over 2}}
-{1\over \sqrt{2}}\beta^{n_0+1+{1\over 2}}\\
&=&
{1\over \sqrt{2}}\left (\alpha^{n_0+1+{1\over 2}}-
\beta^{n_0+1+{1\over 2}}\right )+
{1\over \sqrt{2}}\left (\alpha^{n_0-{1\over 2}}-
\beta^{n_0-{1\over 2}}\right )\\
&=&
{1\over \sqrt{2}}\left (\alpha^{n_0+1+{1\over 2}}-
\beta^{n_0+1+{1\over 2}}\right )+
{1\over \sqrt{2}}\left (\alpha^{n_0-1+{1\over 2}}-
\beta^{n_0-1+{1\over 2}}\right )\\
&=&u_{n_0+1}+u_{n_0-1}.
\end{eqnarray*}
Because by the induction hypothesis $u_{n_0-1}$ is a positive integer, it
follows that $u_{n_0+1}$ is a positive integer. This ends the induction
proof.

\noindent
As a side result we get that $u_n$ obeys the recurrence
$$u_n=4u_{n-1}-u_{n-2},\quad u_0=1,\,u_1=5,$$
since $\alpha^{3\over 2}-\beta^{3\over 2}=\sqrt{50}=5\sqrt{2}$. $u_n$ is
sequence A001834 in \cite{sloane}.

\noindent
It remains to prove that $u_n$ is odd which is the same as to prove that
${L_{2n+1}\over 2}$ is even, which means that $L_{2n+1}$ is a multiple of 4.
Again this will be done by strong induction.
For $n=0$ we have $L_1=4$. Next assume that $L_{2n+1}$ is a multiple of 4
for $n\le n_0$,
that is $L_{2n+1}=4c_n$, where $c_n$ is a positive integer.
Then
\begin{eqnarray*}
L_{2n_0+1}(\alpha^2+\beta^2)&=&4c_{n_0}\cdot 14\\
&=&\alpha^{2n_0+3}+\alpha^{2n_0+1}\beta^2+\alpha^2\beta^{2n_0+1}+
\beta^{2n_0+3}\\
&=&\alpha^{2(n_0+1)+1}+\beta^{2(n_0+1)+1}+\alpha^{2(n_0-1)+1}+
\beta^{2(n_0-1)+1}\\
&=&L_{2(n_0+1)+1} +L_{2(n_0-1)+1}\\
&=&L_{2(n_0+1)+1} +4c_{n_0-1}.
\end{eqnarray*}
Then
$$L_{2(n_0+1)+1}=4c_{n_0}\cdot 14-4c_{n_0-1},$$
that is, a multiple of 4.

\noindent
It follows that
$$b_n={u_n-3\over 2}.$$
Hence
\begin{eqnarray*}
4b_{n-1}-b_{n-2}+3&=&2u_{n-1}-6-{1\over 2}u_{n-2}+{3\over 2}+3\\
&=&{4u_{n-1}-u_{n-2}-3\over 2}\\
&=&{u_n-3\over 2}\\
&=&b_n.
\end{eqnarray*}
The initial values are $b_0=-1,\,b_1=1$.
It follows that $b_n=w(3,\,-1,\,1)$.
The initial values are
$-1,\,1,\,8,\,34,\,131,\,493,\,\ldots$.

\noindent
If we solve Equation~\ref{eq:equazione} for $r$, being $r$ a positive
quantity, we obtain
$$r={\sqrt{3}\sqrt{11+12s^2+4s^2}-3\over 6}.$$
If we define
$$v_n^2=3(11+12b_n+4b_n^2),$$
after insertion of the value of $b_n$ we get
$$v_n^2=3(u_n^2+2),$$
that is
$$v_n^2={3\over 2}(L_{2n+1}+2).$$
Essentially repeating the previous reasoning
we find that $v_n$ is a positive integer $\forall\, n$ and we have
the recurrence
$$v_n=4v_{n-1}-v_{n-2},\quad v_0=3,\,v_1=9,$$
that is $v_n=w(0,\,3,\,9)$. Also we get the closed form
$$v_n=\sqrt{{3\over 2}}\left (\alpha^{n+{1\over 2}}+\beta^{n+{1\over 2}}
\right ).$$
Furthermore
$$v_n\equiv 3\pmod{6}.$$
This can be proved easily by strong induction, noting that $v_0=3$ and
assuming that for all $n\le n_0$ the former equation holds.
Writing $v_{n_0}=6k_1+3,\,v_{n_0-1}=6k_2+3$ we have
\begin{eqnarray*}
v_{n_0+1}&=&4v_{n_0}-v_{n_0-1}\\
&=&4(6k_1+3)-(6k_2+3)\\
&=&6(4k_1-k_2)+9,
\end{eqnarray*}
so that
\begin{eqnarray*}
v_{n_0+1}-3&=&6(4k_1-k_2)+6\\
&=&6(4k_1-k_2+1).
\end{eqnarray*}
And finally we have
$$a_n={v_n-3\over 6}.$$

\noindent
It is interesting the fact that the sequence $u_n$ is a bisection of
sequence A002531 in \cite{sloane} that gives the numerators of continued
fraction convergents to $\sqrt{3}$, a result that as far as we know is new.
Denoting by $z_n$ the sequence of the convergents we are going to prove
$$u_n=z_{2n+1}.$$
From the comments to A002531 in \cite{sloane} we have that
\begin{equation}
z_{2n+1}=2z_{2n}+z_{2n-1},\quad z>0,
\end{equation}
\begin{equation}
z_{2n}={\alpha^n+\beta^n\over 2}.
\end{equation}
Now
\begin{eqnarray*}
u_n-u_{n-1}&=&
{1\over \sqrt{2}}\left (\alpha^{n+{1\over 2}}-\beta^{n+{1\over 2}}\right )
-{1\over \sqrt{2}}\left (\alpha^{n-{1\over 2}}-\beta^{n-{1\over 2}}\right )\\
&=&
{1\over \sqrt{2}}\left (\alpha^n(\alpha^{1\over 2}-\alpha^{-{1\over 2}})+
\beta^n(\beta^{-{1\over 2}}-\beta^{1\over 2})\right )\\
&=&
{1\over \sqrt{2}}\left (\alpha^n(\alpha^{1\over 2}-\beta^{1\over 2})+
\beta^n(\alpha^{1\over 2}-\beta^{1\over 2})\right )\\
&=&
{1\over \sqrt{2}}(\alpha^n+\beta^n)\left (\alpha^{1\over 2}-\beta^{1\over 2}
\right ).
\end{eqnarray*}
Now from
$$\left (\alpha^{1\over 2}+\beta^{1\over 2}\right )^2=\alpha+\beta+2=6$$
we get
$$\alpha^{1\over 2}+\beta^{1\over 2}=\sqrt{6}.$$
Also from
$$\left (\alpha^{1\over 2}+\beta^{1\over 2}\right )
\left (\alpha^{1\over 2}-\beta^{1\over 2}\right )=\alpha-\beta=2\sqrt{3}$$
we get
$$\alpha^{1\over 2}-\beta^{1\over 2}=\sqrt{2},$$
so that
$$u_n-u_{n-1}=L_n=2z_{2n}=z_{2n+1}-z_{2n-1}.$$
So we can write
$$u_n-z_{2n+1}=u_{n-1}-z_{2n-1}.$$
Setting
$$d_n=u_n-z_{2n+1}$$
we have the recurrence
$$d_n=d_{n-1}$$
with $d_0=u_0-z_1=1-1=0$. Hence $d_n=0$ and $u_n=z_{2n+1}$.

\noindent
In a similar way we obtain
$$v_n-v_{n-1}=6F_n,$$
where $F_n=w(0,\,0,\,1)$, the analogous of the Fibonacci numbers for this
type of recurrences. $F_n$ is
sequence A001353 in \cite{sloane}.

\end{document}